# Efficient partially replicated block designs with each replication number one or two

R.A. Bailey  and  Rahul Mukerjee*
School of Mathematics and Statistics     Indian Institute of Management Calcutta
University of St Andrews                 Joka, Diamond Harbour Road
St Andrews, KY16 9SS, Fife, UK           Kolkata 700 104, India
E-mail: rab24@st-andrews.ac.uk           Email: rmuk0902@gmail.com
                                         * Corresponding author

*Abstract*: We investigate block designs, under the *A*- and *MV*-criteria, when each treatment can have only one or two replications due to resource constraints, as can happen, for example, in early generation varietal trials. While these are commonly known as partially replicated designs, a key new feature of the present work is that no restriction about a constant block size is imposed on the subdesign consisting of the twice replicated treatments. This makes the derivation more challenging but allows us to entertain a wider class of competing designs and hence increases the flexibility of the results. Considering all treatments as equally important, design-independent, sharp lower bounds on the *A*- and *MV*-criteria are derived, so as to find highly efficient designs over this wider class. The roles of (a) linked block designs, (b) designs in an online catalog `designtheory.org`, and (c) partially balanced incomplete block designs, or duals thereof, as adapted to our setup, are explored at length. Illustrative examples are presented.

## 1. Introduction

We consider block designs when the available resources allow very few replications for the treatments, permitting two replications for a relatively small number of them and only one replication for the remaining vast majority of them. The focus is on obtaining highly efficient designs for inference on pairwise treatment contrasts, under the *A*- and *MV*-criteria and over a wide class of competing designs. Recall that the *A*-criterion aims at estimating these contrasts with the minimum possible average variance, while the *MV*-criterion is about minimizing the maximum of these variances.

An experimental situation as above can arise, for example, in agriculture when interest lies in comparing a large number of early generation varieties which represent the treatments. Available resources for these, like the quantity of seed, are often quite limited, permitting the use of two plots for only some of them and a single plot for the rest; see e.g., Haines (2021) who called these partially replicated designs and cited further references including Williams et al. (2011, 2014), among others. As elaborated in the next section, however, a major difference between our approach and that in Haines (2021) is that, unlike in some other earlier work including hers, we do not assume a constant block size for the subdesign formed by the twice replicated treatments. This adds considerably to the complexity of our derivation but leads to results in the framework of a wider class of designs, thus allowing for more flexibility in the choice of efficient designs.

It is worthwhile to mention here that the term partially replicated or ($p$-rep) design was first used by Cullis et al. (2006), but their setup and approach are different from ours. They considered row-column designs along with random genotype effects, whereas we study block designs under a fixed effects model.

Bailey and Cameron (2013) and Bailey and Sajjad (2021) investigated a graph theoretic approach for finding optimal block designs that keep all treatment contrasts estimable using a minimal or nearly minimal number of plots. This is similar to our situation in spirit, but a key difference is that their optimal designs can have appreciably more than two replications for some treatments, even to the extent of having a treatment that appears in all blocks. So, their results, as they stand, do not apply to our setup where all treatments have severely constrained replications, with no more than two replications permitted for any of them. Interestingly, however, our results show that the duals of some of their designs can perform very well in the present context.

Section 2 describes our setup in more detail. It also presents some notation and a preliminary result. Then, we obtain an expression for the *A*-criterion in Section 3. A sharp lower bound on the *A*-criterion is derived in Section 4. It is employed in this and the next two sections to explore designs with high *A*-efficiency at length. Section 7 dwells on highly efficient designs under the *MV*-criterion. Section 8 shows that our designs remain highly efficient also under two variants of the *A*-criterion. In Section 9, we examine how well our approach works for practically relevant large designs. A guideline on the use of our results is also presented there. The paper is concluded in Section 10 with an indication of several open issues that arise from the present work. Finally, some further details appear in the appendices.

**2. Setup, notation and a proposition**

Consider a block design $d_0$ involving $b$ blocks, each of size $k$, and $u + w$ treatments of which $u$ are replicated twice each and $w$ are replicated once each. Let $U$ and $W$ denote the sets of these $u$ and $w$ treatments, respectively. The numbers $u$, $w$, $b$ and $k$ are pre-specified by the experimental requirements. Clearly, they satisfy $bk = 2u + w$. There are $v = u + w$ treatments altogether and, as mentioned in the introduction, we focus on the situation where $w$ is much larger than $u$.

Interest lies in inference on all the $v(v-1)/2$ pairwise treatment contrasts. Hence only connected designs are considered. Let $d$ be the subdesign of $d_0$ that contains the treatments in $U$. Thus, each treatment is replicated twice in $d$. Clearly, $d_0$ is connected if and only if $d$ is connected and spread out over all the $b$ blocks. Then the estimation space of $d$ has rank $u + b - 1$, while $d$ involves $2u$ plots. So, connectedness necessitates $2u \geq u + b - 1$, i.e., $u \geq b - 1$, which is, hereafter, assumed to hold. The subdesign $d$ determines $d_0$ uniquely up to renaming of the treatments in $W$, and plays a crucial role in the sequel. Note also that $d_0$ is binary if and only if so is $d$. However, unless otherwise specified, no assumption about binarity is made. For clarity, the design $d_0$, involving all the $v$ $(= u + w)$ treatments, will occasionally be referred to as the full design.

For $j = 1,\ldots, b$, suppose the $j$th block of $d_0$ has $s_j$ and $k - s_j$ plots that are populated by the treatments in $W$ and $U$, respectively. Then $d$ has block sizes $k - s_j$, $j = 1,\ldots, b$, and $0 \leq s_j \leq k - 1$ for every $j$, as connectedness requires $d$ to be spread out over all the $b$ blocks. Moreover,

$$w = s_1 + \cdots + s_b. \tag{1}$$

Along the lines of Haines (2021, p. 415), the treatments from $U$ are allocated to plots randomly within each block of $d_0$, and then those from $W$ are allocated randomly to the remaining plots. We consider a fixed effects additive linear model, where the yields from the $bk$ plots of $d_0$ are assumed to be uncorrelated, with a common variance $\sigma^2$. Write $\tau_1,\ldots, \tau_u$ for the fixed effects of the $u$ treatments in $U$. Also, for each $j$ with $s_j \geq 1$, let $\rho_{jh}$ be the fixed effect of the $h$th of the $s_j$ treatments in $W$ that appear in the $j$th block, $h = 1,\ldots, s_j$. By (1), each $\rho_{jh}$ corresponds to one treatment in $W$.

Let $N = (n_{ij})$ be the $u \times b$ incidence matrix of the subdesign $d$, where $n_{ij}$ is the number of times the $i$th treatment in $U$ appears in the $j$th block. Then the dual, $\tilde{d}$, of $d$ has incidence matrix $N^T$, of order $b \times u$, where the superscript $T$ denotes transpose. Note that $\tilde{d}$ has a constant block size two as $d$ is two-replicate. Since $d$ has block sizes $k - s_j$, $j = 1,\ldots, b$, writing $I_u$ and $I_b$ for the identity matrices of orders $u$ and $b$, and $S = \text{diag}(s_1,\ldots, s_b)$, the intrablock information matrices of $d$ and $\tilde{d}$ are then given, respectively, by

$$C = 2I_u - N(kI_b - S)^{-1}N^T, \qquad \tilde{C} = kI_b - S - (1/2)N^TN. \tag{2}$$

Let $C^+$ and $\tilde{C}^+$ denote their respective Moore-Penrose inverses (Rao, 1973, p. 26). Write $e_1,\ldots, e_u$ for the unit vectors of order $u$, and $\tilde{e}_1,\ldots, \tilde{e}_b$ for the unit vectors of order $b$. Also, let

$$\xi_{ij} = \tilde{e}_j - (1/2)N^Te_i, \qquad i = 1,\ldots, u, \qquad j = 1,\ldots, b. \tag{3}$$

We are now in a position to present the following proposition, which shows the variances of best linear unbiased estimators (BLUEs) of pairwise treatment contrasts, i.e., pairwise contrasts involving $\tau_1,\ldots, \tau_u$ and the $\rho_{jh}$, and will be the starting point of our results under the $A$- and $MV$-criteria.

**Proposition 1**. *The following hold for the full design $d_0$:*

(a) *For $1 \leq i < i^* \leq u$, the BLUE of $\tau_i - \tau_{i^*}$ has variance $\sigma^2(e_i - e_{i^*})^T C^+ (e_i - e_{i^*})$.*

(b1) *For every j, with $s_j \geq 2$, and $1 \leq h < h^* \leq s_j$, the BLUE of $\rho_{jh} - \rho_{jh^*}$ has variance $2\sigma^2$.*

(b2) *For every $j, j^*$, with $j < j^*$ and $s_j, s_{j^*} \geq 1$, and $1 \leq h \leq s_j$, $1 \leq h^* \leq s_{j^*}$, the BLUE of $\rho_{jh} - \rho_{j^*h^*}$ has variance $\sigma^2 \{2 + (\tilde{e}_j - \tilde{e}_{j^*})^T \tilde{C}^+ (\tilde{e}_j - \tilde{e}_{j^*})\}$.*

(c) *For every j, with $s_j \geq 1$, $1 \leq h \leq s_j$ and $1 \leq i \leq u$, the BLUE of $\tau_i - \rho_{jh}$ has variance $\sigma^2 \{(3/2) + \xi_{ij}^T \tilde{C}^+ \xi_{ij}\}$.*

Observe that parts (a)-(c) of Proposition 1 collectively take care of all the $v(v-1)/2$ pairwise treatment contrasts. Since $v = u + w$, this follows because (a) accounts for $u(u-1)/2$ of these contrasts, while (b1) and (b2) together take care of $\sum_{j=1}^{b} s_j(s_j - 1)/2 + \Sigma\Sigma_{1 \leq j < j^* \leq b} s_j s_{j^*} = w(w-1)/2$ of them by (1), and (c) covers $\sum_{i=1}^{u} \sum_{j=1}^{b} s_j = uw$ of them.

Proposition 1 can be obtained from Mukerjee (2025, Proposition 2), parts of which also appear in Herzberg and Jarrett (2007) and Haines (2021). Compared to these papers, however, the present design problem has two major differences which make the mathematical development here quite different:

(i) As mentioned in the introduction, unlike in these papers, the subdesign $d$ for the treatments in $U$ is not assumed to have a constant block size. The block sizes of $d$ are $k - s_j$, $j = 1,\ldots, b$, where $s_1,\ldots, s_b$ are not pre-fixed and choosing these design variables is an integral part of our design problem. Even if the highly efficient designs obtained later turn out to have equal or nearly equal $s_1,\ldots, s_b$, such efficiency has to be established in the entire class of designs for a given $(u, w, b, k)$, without any assumption on $s_1,\ldots, s_b$. Thus, our results apply to a wider class of designs compared to the earlier papers.

(ii) All the $v$ ($= u + w$) treatments are equally important in our setup, which justifies use of overall $A$- or $MV$-criteria, as against separate $A$- or $MV$-criteria for various kinds of treatment comparisons in Herzberg and Jarret (2007) and Mukerjee (2025). Moreover, we shall obtain, for the first time in the present setup, design-independent lower bounds on the overall $A$- or $MV$-criteria, leading to efficiency measures which facilitate identification of good designs in an objective manner.

**3. Expression for the $A$-criterion**

Theorem 1 below combines the variances in Proposition 1 (a)-(c) and leads to an expression for the overall $A$-criterion. Although Proposition 1 appears in some earlier work cited in the last section, Theorem 1 is new. We re-emphasize that, in contrast to such earlier results, Theorem 1 involves the design variables $s_1,\ldots, s_b$ and this happens via not only $\text{tr}(\tilde{C}^+ S)$ but also $\text{tr}(C^+)$ and $\text{tr}(\tilde{C}^+)$; see (2).

**Theorem 1**. *The variances of the BLUEs of all the $v(v-1)/2$ pairwise treatment contrasts in the full design $d_0$ sum to $\sigma^2 A(d_0)$, where*

$$A(d_0) = (1/2)w(3u + 2w - b - 1) + u\mathrm{tr}(C^+) + (1/2)k\{w\mathrm{tr}(\tilde{C}^+) + b\mathrm{tr}(\tilde{C}^+ S)\}.$$

*Proof.* Several identities will be useful in the proof. First note that,

$$\sum_{j=1}^{b} s_j \tilde{e}_j = s, \qquad \sum_{j=1}^{b} s_j \tilde{e}_j \tilde{e}_j^T = S, \tag{4}$$

where $s = (s_1, \ldots, s_b)^T$ and, as before, $S = \mathrm{diag}(s_1, \ldots, s_b)$. Next, write $1_u$ and $1_b$, respectively, for the $u \times 1$ and $b \times 1$ vectors of ones, and observe that

$$\sum_{i=1}^{u} N^T e_i = N^T 1_u = k 1_b - s, \qquad \sum_{i=1}^{u} N^T e_i e_i^T N = N^T N = 2(kI_b - S - \tilde{C}). \tag{5}$$

In the first identity of (5), we recall that the subdesign $d$, with incidence matrix $N$, has block sizes $k - s_j$, $j = 1, \ldots, b$, while we use (2) in the second identity. Furthermore,

$$C^+ 1_u = 0, \qquad \tilde{C}^+ 1_b = 0, \qquad \tilde{C}^+ \tilde{C} = I_b - b^{-1} J_{bb}, \tag{6}$$

where $J_{bb}$ is the $b \times b$ matrix of ones. These hold because $C^+$ and $\tilde{C}^+$ are Moore-Penrose inverses and the subdesign $d$ is connected.

By (6), the contribution of the $u(u-1)/2$ variances in Proposition 1(a) to $A(d_0)$ is

$$A_{UU} = \sum\sum_{1 \le i < i* \le u} (e_i - e_{i*})^T C^+ (e_i - e_{i*}) = u\mathrm{tr}(C^+),$$

invoking a well-known result on the average variance of BLUEs of elementary treatment contrasts in a connected block design; see, e.g., Williams and Piepho (2015). Next, by (1) and (4), the contribution of the $w(w-1)/2$ variances in Proposition 1(b1, b2) to $A(d_0)$ equals

$$A_{WW} = \sum_{j=1}^{b} s_j(s_j - 1) + \sum\sum_{1 \le j < j* \le b} s_j s_{j*} \{2 + (\tilde{e}_j - \tilde{e}_{j*})^T \tilde{C}^+ (\tilde{e}_j - \tilde{e}_{j*})\}$$

$$= w(w-1) + \sum\sum_{1 \le j < j* \le b} s_j s_{j*} (\tilde{e}_j - \tilde{e}_{j*})^T \tilde{C}^+ (\tilde{e}_j - \tilde{e}_{j*})$$

$$= w(w-1) + (1/2)\mathrm{tr}\{\tilde{C}^+ \sum_{j=1}^{b} \sum_{j*=1}^{b} s_j s_{j*} (\tilde{e}_j - \tilde{e}_{j*})(\tilde{e}_j - \tilde{e}_{j*})^T\}$$

$$= w(w-1) + \mathrm{tr}\{\tilde{C}^+ (wS - ss^T)\}.$$

Finally, by (1) and (3)-(6), the contribution of the $uw$ variances in Proposition 1(c) to $A(d_0)$ is given by

$$A_{UW} = \sum_{i=1}^{u} \sum_{j=1}^{b} s_j \{(3/2) + \xi_{ij}^T \tilde{C}^+ \xi_{ij}\} = (3/2)uw + \mathrm{tr}(\tilde{C}^+ \sum_{i=1}^{u} \sum_{j=1}^{b} s_j \xi_{ij} \xi_{ij}^T)$$

$$= (3/2)uw + \mathrm{tr}[\tilde{C}^+ \sum_{i=1}^{u} \sum_{j=1}^{b} s_j \{\tilde{e}_j - (1/2)N^T e_i\}\{\tilde{e}_j - (1/2)N^T e_i\}^T]$$

$$= (3/2)uw + \mathrm{tr}[\tilde{C}^+ \{uS - (1/2)(k1_b - s)s^T - (1/2)s(k1_b - s)^T + (1/2)w(kI_b - S - \tilde{C})\}]$$

$$= (3/2)uw + \text{tr}[\tilde{C}^+\{(u-w/2)S + ss^T + (1/2)w(kI_b - \tilde{C})\}]$$

$$= (3/2)uw + \text{tr}[\tilde{C}^+\{(u-w/2)S + ss^T + (1/2)wkI_b\}] - (1/2)w\text{tr}(I_b - b^{-1}J_{bb})$$

$$= (3/2)uw - (1/2)w(b-1) + \text{tr}[\tilde{C}^+\{(u-w/2)S + ss^T + (1/2)wkI_b\}].$$

Summing $A_{UU}$, $A_{WW}$ and $A_{UW}$, the result follows if we recall that $bk = 2u + w$. □

We remark that, analogously to $A_{UU}$ in the above proof, we have $A(d_0) = v\text{tr}(C_0^+)$, where $C_0^+$ is the Moore-Penrose inverse of the intrablock information matrix $C_0$ of the full design $d_0$. Similarly, following Cullis et al. (2006, p. 385),

$$A_{WW} = w\text{tr}(C_{0w}^+) - 1_w^T C_{0w}^+ 1_w,$$

where $C_{0w}^+$ is the principal submatrix of $C_0^+$ corresponding to the $w$ treatments in $W$ and $1_w$ is the $w\times 1$ vector of ones. While these alternative expressions for $A(d_0)$ and $A_{WW}$ are apparently simpler than their counterparts in Theorem 1, the latter are advantageous for our purpose and hence used in the sequel because they allow us to exploit the properties of the subdesign $d$ which is much smaller than $d_0$.

By Theorem 1, the A-criterion is given by $2A(d_0)/\{v(v-1)\}$, ignoring the multiplicative constant $\sigma^2$, and it suffices to work with $A(d_0)$ for exploring this criterion. This theorem enables us to compare connected full designs like $d_0$, which involve $b$ blocks, each of size $k$, along with $u$ treatments each replicated twice and $w$ treatments each replicated once. Given $u$, $b$ and $k$, which uniquely determine $w = bk - 2u$, we denote this class of competing designs by $D_0 \equiv D_0(u, b, k)$. An A-optimal design in $D_0$ is one that minimizes the average variance or equivalently, $A(d_0)$ in Theorem 1. However, the presence of $s_1,\ldots, s_b$ in $A(d_0)$ complicates the exact A-optimality problem in general. For instance, although intuition suggests taking these design variables as equal or nearly equal, theoretical justification for doing so remains elusive. Hence, we now focus on designs that are highly A-efficient. This calls for obtaining a sharp lower bound on $A(d_0)$ over $D_0$. The following lemma will be useful for this purpose and also later in the paper.

**Lemma 1**. (a) *For any nonnegative definite (nnd) matrix $\Delta$ of order $b$,*

$$\text{tr}(\tilde{C}\Delta)\text{tr}(\tilde{C}^+\Delta) \geq \{\text{tr}(\Delta) - b^{-1}1_b^T \Delta 1_b\}^2.$$

(b) *Furthermore, if $\Delta = \text{diag}(\delta_1, \ldots, \delta_b)$, where $\delta_1, \ldots, \delta_b$ are nonnegative, then*

$$\text{tr}(\tilde{C}\Delta) \leq (1/2)\Sigma_{j=1}^b (k - s_j)\delta_j.$$

*Proof.* (a) Since $d$ is connected and $\tilde{C}^+$ is the Moore-Penrose inverse of $\tilde{C}$, spectral decompositions yield

$$\text{tr}(\tilde{C}\Delta) = \Sigma_{j=1}^{b-1}\theta_j q_j^T \Delta q_j, \qquad \text{tr}(\tilde{C}^+\Delta) = \Sigma_{j=1}^{b-1}\theta_j^{-1} q_j^T \Delta q_j,$$

where $\theta_1,\ldots,\theta_{b-1}$ are the positive eigenvalues of $\tilde{C}$, with $q_1,\ldots,q_{b-1}$ as the corresponding orthonormal eigenvectors which satisfy

$$\Sigma_{j=1}^{b-1} q_j q_j^T = I_b - b^{-1} 1_b 1_b^T.$$

Hence (a) follows by Cauchy-Schwarz inequality.

(b) For $j = 1,\ldots, b$, by (2), the $j$th diagonal element of $\tilde{C}$ is $k - s_j - (1/2)\Sigma_{i=1}^{u} n_{ij}^2$, which does not exceed $(1/2)(k - s_j)$, because $n_{ij}^2 \geq n_{ij}$ and $\Sigma_{i=1}^{u} n_{ij} = k - s_j$. Hence (b) follows. □

## 4. Lower bound on the $A$-criterion and $A$-efficient designs

For every design $d_0$ in $D_0$, we have

$$\text{tr}(C^+) \geq (u-1)^2/\text{tr}(C) \geq (u-1)^2/(2u-b), \quad \text{tr}(\tilde{C}^+) \geq (b-1)^2/\text{tr}(\tilde{C}) \geq (b-1)^2/u, \tag{7}$$

$$\text{tr}(\tilde{C}^+ S) \geq w(b-1)^2/(bu). \tag{8}$$

The inequalities in (7) follow readily from (2) and the arithmetic mean harmonic mean inequality, as $n_{ij}^2 \geq n_{ij}$. To see the truth of (8), note that by Lemma 1(b),

$$\text{tr}(\tilde{C}S) \leq (1/2)\Sigma_{j=1}^{b}(k - s_j)s_j \leq wu/b,$$

because $bk = 2u + w$, and by (1),

$$\Sigma_{j=1}^{b}(k - s_j)s_j = kw - \Sigma_{j=1}^{b} s_j^2 \leq kw - (w^2/b) = (w/b)(bk - w).$$

Also, by (1), $\text{tr}(S) - b^{-1} 1_b^T S 1_b = w(b-1)/b$. Hence, taking $\Delta = S$ in Lemma 1(a), we obtain (8).

By (7) and (8), a design-independent lower bound on $A(d_0)$ in Theorem 1, over the class $D_0$, is obtained as

$$A_{\text{bound}} = (1/2)w(3u + 2w - b - 1) + u\{(u-1)^2/(2u-b)\} + kw\{(b-1)^2/u\}. \tag{9}$$

This sets a lower bound on the $A$-criterion as $2A_{\text{bound}}/\{v(v-1)\}$, and gives a measure of the $A$-efficiency of any design $d_0$ in $D_0$ as

$$A_{\text{eff}}(d_0) = A_{\text{bound}}/A(d_0). \tag{10}$$

The $A$-efficiency measure in (10) is conservative since one or more of the inequalities in (7) and (8) may remain strict over $D_0$ and thus the lower bound in (9) may not be attainable. Nevertheless, this bound is quite sharp in many situations, enabling us to find designs with high $A$-efficiencies. We begin showing this in Subsection 4.1 below via consideration of linked block and nearly linked block designs. Furthermore, as the proportion $f = w/v$ of the singly replicated treatments in $d_0$ is intended to be large, the values of $f$ are also indicated in our examples and tables. By and large, the designs obtained here will be seen to be highly efficient, irrespective of the value of $f$, barring the unrealistic cases where $f$ is too small.

Given any subdesign $d$, we shall write $k_{\max}(d)$ for the largest block size in $d$ and $k_0$ [$\equiv k_0(d)$] for the smallest possible value of the block size, $k$, of the corresponding full design $d_0$, so as to ensure $w > 0$. Clearly, $k_0$ equals $k_{\max}(d) + 1$ if $d$ has a constant block size, and $k_{\max}(d)$ otherwise.

*4.1. Use of linked block and nearly linked block designs*

**Example 1** (linked block design). To show an immediate application of (10), let $u = \lambda\binom{b}{2}$, $\lambda \geq 1$, $b \geq 3$, and suppose the subdesign $d$ is chosen as a two-replicate linked block design, i.e., its dual $\tilde{d}$ is a balanced incomplete block (BIB) design with $\binom{b}{2}$ distinct blocks $\{j, j^*\}$, $1 \leq j < j^* \leq b$, each repeated $\lambda$ times. Then $d$ has a constant block size $\lambda(b-1)$, $k_0 = \lambda(b-1) + 1$, $k \geq k_0$, and $s_1 = \cdots = s_b = k - \lambda(b-1)$.

In this situation, $\tilde{C} = (\lambda b/2)(I_b - b^{-1}J_{bb})$ and equality holds both in (8) and the second inequality of (7). Furthermore, the positive eigenvalues of $C$ are $b/(b-1)$ and 2, with multiplicities $b - 1$ and $u - b$, respectively (cf. Roy, 1958). Therefore, $\text{tr}(C^+) = \{2(b-1)^2 + b(u-b)\}/(2b)$, and one can check that equality holds also in the first inequality of (7) if and only if $u = b$ (= 3), in which case the subdesign $d$ leads to an A-optimal full design $d_0$. Then $d$ itself is a BIB design and, in fact, the only BIB design which replicates each treatment twice. Generally, with $d$ chosen as a linked block design as above, even when $\text{tr}(C^+)$ exceeds the lower bound on it in (7), the ratio of the two is close to 1, so that $d_0$ has high A-efficiency. Indeed, for $\lambda = 1, 2, 3$, and $b = 3, \ldots, 8$, computations using (10) show that $A_{\text{eff}}(d_0) \geq 0.97$ for $k \geq k_0$, and $A_{\text{eff}}(d_0)$ exceeds 0.99 quickly with increase in $k$, e.g., if $\lambda = 2$ and $b = 4$, then $A_{\text{eff}}(d_0)$ equals 0.986 for $k = k_0 = 7$, and exceeds 0.99 for $k \geq 8$. The same pattern persists for $b \geq 9$ as well.

Here, for instance, with $\lambda = 1$, if $k = 15$, then the proportion, $f$, of singly replicated treatments equals 0.85, 0.80 and 0.75 for $b = 5$, 6 and 7, respectively. The picture is even better for $k = 20$, when the corresponding values of $f$ are 0.89, 0.86 and 0.82. Also, with $\lambda = 1$, for $b = 5, 6, 7$, we have $v = 65, 75, 84$ if $k = 15$, and $v = 90, 105, 119$ if $k = 20$. So, both values of $k$ are relatively small compared to $v$, as they should be. □

In the spirit of Example 1, if $u$ is nearly equal to $\lambda\binom{b}{2}$ for some $\lambda$, then one can anticipate that a nearly linked block design should perform very well. The next example corroborates this point.

**Example 2**. (a) With $\lambda = 1$ and $b = 6$, a two-replicate linked block design for 15 treatments and 6 blocks is given by the blocks

$$\{1, 2, 3, 4, 5\}, \{1, 6, 7, 8, 9\}, \{2, 6, 10, 11, 12\},$$
$$\{3, 7, 10, 13, 14\}, \{4, 8, 11, 13, 15\}, \{5, 9, 12, 14, 15\}.$$

Now, let $(u, b) = (13, 6)$, and consider a subdesign $d$ obtained by deleting the two treatments 1 and 15 from the above linked block design. The block sizes of $d$ differ among themselves by at most unity and hence the same holds for $s_1,\ldots, s_6$. From (10), one can check that the resulting $A_{\text{eff}}(d_0)$ equals 0.952, 0.971 and 0.978 for $k = 5, 6$ and 7, respectively, and exceeds 0.98 for $k \geq 8$.

(b) Let $(u, b) = (17, 6)$. Consider a subdesign $d$ obtained by adding a treatment 16 to the first two blocks and another treatment 17 to the last two blocks of the linked block design in (a). Then $A_{\text{eff}}(d_0)$ equals 0.952 and 0.974 for $k = 6$ and 7, respectively, and is greater than 0.98 for $k \geq 8$.

(c) Let $(u, b) = (8, 4)$. Suppose the subdesign $d$ consists of the blocks $\{1, 2, 3, 7\}$, $\{1, 4, 5, 7\}$, $\{2, 4, 6, 8\}$, $\{3, 5, 6, 8\}$, i.e., $d$ is obtained by adding treatments 7 and 8 to a linked block design. Then $A_{\text{eff}}(d_0)$ equals 0.975 for $k = 5$ and exceeds 0.98 for $k \geq 6$.

In parts (a)-(c) above, $f$ equals 0.83, 0.77, 0.85 for $k = 15$, and 0.88, 0.84, 0.89 for $k = 20$. □

**Remark 1**. The subdesign for $u = b = 3$ in Example 1 ensures $A_{\text{eff}}(d_0) = 1$, for each $k \geq k_0 (= 3)$. For every other subdesign $d$ in Examples 1 and 2 and also in the rest of this paper, we find that $A_{\text{eff}}(d_0)$ increases with $k$. This is reassuring as, given $d$, a larger $k$ is associated with a larger $w$ which is typically of greater practical interest. Moreover, this facilitates the presentation of our findings. To that end, for any given $d$ and $\alpha = 0.90, 0.95$ and 0.98, let $k_\alpha$ denote the smallest integer, if any, such that $k_0 \leq k_\alpha \leq 20$ and $A_{\text{eff}}(d_0) \geq \alpha$, whenever $k \geq k_\alpha$, e.g., $(k_0, k_{0.90}, k_{0.95}, k_{0.98}) = (5, 5, 5, 8)$ in Example 2(a). If no such $k_\alpha$ exists, then we write $k_\alpha = \text{X}$, to signify that $A_{\text{eff}}(d_0)$ does not grow adequately fast with $k$ so as to reach $\alpha$. The $k_\alpha$ can be readily obtained using the expression for $A_{\text{eff}}(d_0)$ in (10). □

The ideas in Examples 1 and 2, based on linked block designs or modifications thereof, work when $u$ is equal or close to $\lambda\binom{b}{2}$ for some $\lambda$. In the next two sections, we examine other helpful approaches for choosing the subdesign $d$.

## 5. Use of an online design resource

The website `designtheory.org` is an important design resource which is worth exploration. The catalog `v-b-k-phi` there shows $A$-optimal block designs in the binary class and is of particular relevance to us. Each of these designs has a constant block size and is indexed by an ordered triplet; for instance, 4-7-2 represents a design for 4 treatments in 7 blocks each of size 2. It is not hard to see that among these designs, the ones with block size two are $A$-optimal even when non-binary designs are entertained. Moreover, their duals are two-replicate and hence qualify as the subdesign $d$ in our setup. The $d$ so obtained for any $u$ and $b$ is denoted by $d_{\text{cat}}(u, b)$, i.e., $d_{\text{cat}}(u, b)$ is the dual of the $A$-optimal block design $b$-$u$-2 in the catalog v-b-k-phi. The catalog extends over the range $u + b \leq 16$ in our notation and from a practical viewpoint, we consider $u, b \geq 4$. For every such $u, b$, the catalog reports a unique A-optimal b-u-2 and hence there is no ambiguity about $d_{\text{cat}}(u, b)$. Some of these subdesigns are linked block or nearly linked block designs, e.g., $d_{\text{cat}}(8, 4)$ is the same as the $d$ in Example 2(c), but many are not. Generally, the $d_{\text{cat}}(u, b)$ hold promise as they minimize $\text{tr}(\widetilde{C}^+)$, their duals being A-optimal block designs. At the same time, $A(d_0)$ in Theorem 1, and hence our A-criterion, are much more complex because they involve all three of $\text{tr}(C^+)$, $\text{tr}(\widetilde{C}^+)$ and $\text{tr}(\widetilde{C}^+ S)$. Consequently, a closer examination of the $d_{\text{cat}}(u, b)$ is warranted.

For each $d_{\text{cat}}(u, b)$, over $u, b \geq 4$, $u + b \leq 16$, we first employ (10) to study the A-efficiency, $A_{\text{eff}}(d_0)$, of the full design $d_0$ and obtain $k_0$, $k_{0.90}$, $k_{0.95}$ and $k_{0.98}$ introduced in the last section. These quantities, shown in Table 1 below, reveal that for many $(u, b)$, $A_{\text{eff}}(d_0)$ quickly reaches 0.95 or 0.98 with increase in $k$, but at the same time, there are quite a few other $(u, b)$ where $A_{\text{eff}}(d_0)$ falls short of 0.95 or even 0.90, at $k = 20$, as signified by $k_{0.95} = $ X or $k_{0.90} = $ X. In all these latter cases, $d_{\text{cat}}(u, b)$ has average block size $2u/b < 3$, and even has several blocks of size only one when $u = b - 1$, as happens when it is the dual of a queen-bee design (see e.g., Bailey and Cameron, 2013, Subsection 9.1). Thus, the apparently unimpressive performance of $d_{\text{cat}}(u, b)$ in these situations can be either (i) an indicator of the reality, or (ii) due to possible over-conservativeness of $A_{\text{eff}}(d_0)$ itself for such small block sizes of the subdesign.

Computations indicate the truth of the second possibility above. For every $(u, b)$ in Table 1, an enumeration of all possible binary subdesigns via consideration of their duals shows that, among them, $d_{\text{cat}}(u, b)$ always minimizes $A(d_0)$, for block size $k = k_0, 5, 10, 15, 20$, or $k = k_0, 10, 15, 20$ if $k_0 > 5$. The same holds also for numerous other values of $k$ that we considered, without ever coming across any better subdesign. Thus, in reality, throughout Table 1 including the apparently doubtful cases there, $d_{\text{cat}}(u, b)$ has the best performance in the binary competing class under the A-criterion, and hence is expected to work well even if non-binary competitors are entertained. In fact, if $u = b - 1$, then it is not hard to see that all connected subdesigns have to be binary, which further strengthens our enumerative findings. We again stress that, due to the complex form of $A(d_0)$, the facts just noted are far from obvious from mere A-optimality of the dual of any $d_{\text{cat}}(u, b)$ as a block design.

To summarize, we observe that the use of the online catalog v-b-k-phi completely settles the task of finding A-efficient designs for $u, b \geq 4$, $u + b \leq 16$. For any such $(u, b)$, one only has to find the A-optimal block design b-u-2 from v-b-k-phi, consider its dual $d_{\text{cat}}(u, b)$ as the subdesign, and remain assured of the A-optimality of the resulting full design $d_0$ among binary competing designs in general, and in the entire competing class when $u = b - 1$. Moreover, in many situations, $d_{cat}(u, b)$ ensures high values of $A_{eff}(d_0)$, that pertain to all competing designs including the non-binary ones. We further note that throughout Table 1, the proportion, $f$, of singly replicated treatments is 0.75 or higher for $k = 15$, and 0.82 or higher for $k = 20$. In fact, even for $k = 15$, we have $f \geq 0.80$ in 27 out of 29 cases in Table 1.

Table 1. *Values of $k_0, k_{0.90}, k_{0.95}$ and $k_{0.98}$ for subdesigns $d_{cat}(u, b)$*

| $(u, b)$ | $k_0$ | $k_{0.90}$ | $k_{0.95}$ | $k_{0.98}$ | $(u, b)$ | $k_0$ | $k_{0.90}$ | $k_{0.95}$ | $k_{0.98}$ | $(u, b)$ | $k_0$ | $k_{0.90}$ | $k_{0.95}$ | $k_{0.98}$ |
|---|---|---|---|---|---|---|---|---|---|---|---|---|---|---|
| (4, 4) | 3 | 3 | 4 | X | (5, 5) | 3 | 7 | X | X | (8, 6) | 3 | 5 | X | X |
| (5, 4) | 3 | 3 | 6 | X | (6, 5) | 3 | 4 | X | X | (9, 6) | 4 | 4 | 5 | X |
| (6, 4) | 4 | 4 | 4 | 4 | (7, 5) | 3 | 4 | 7 | X | (10, 6) | 4 | 4 | 6 | X |
| (7, 4) | 4 | 4 | 4 | 6 | (8, 5) | 4 | 4 | 5 | X | (6, 7) | 6 | X | X | X |
| (8, 4) | 5 | 5 | 5 | 6 | (9, 5) | 4 | 4 | 5 | 9 | (7, 7) | 3 | X | X | X |
| (9, 4) | 5 | 5 | 5 | 7 | (10, 5) | 5 | 5 | 5 | 5 | (8, 7) | 3 | X | X | X |
| (10, 4) | 6 | 6 | 6 | 6 | (11, 5) | 5 | 5 | 5 | 7 | (9, 7) | 3 | 11 | X | X |
| (11, 4) | 6 | 6 | 6 | 7 | (5, 6) | 5 | X | X | X | (7, 8) | 7 | X | X | X |
| (12, 4) | 7 | 7 | 7 | 7 | (6, 6) | 3 | X | X | X | (8, 8) | 3 | X | X | X |
| (4, 5) | 4 | X | X | X | (7, 6) | 3 | 9 | X | X | | | | | |

Some further details on A-efficient designs, related to this and the previous section, also covering the connection with Bailey and Cameron (2013), are presented in Appendix 1.

**6. Use of PBIB designs and their duals**

The last section discussed at length how the online catalog `v-b-k-phi` facilitates efficient choice of the subdesign $d$, when $u + b \leq 16$. Turning to $u + b > 16$, we now explore the role of partially balanced incomplete block (PBIB) designs in this regard. Many of these designs are promising because of their high A-efficiencies as block designs. Of relevance to us are the PBIB designs which are either

(a) two-replicate and hence can themselves be used as $d$, or

(b) have block size two and hence allow the use of their duals as $d$.

We consider $u + b > 16$ and, for ease in presentation, confine our attention to $u, b \geq 4$, $u + b \leq 30$, excluding the cases where $u = \lambda\binom{b}{2}$, since these have already been covered in Example 1. Over this range of $(u, b)$, we study all subdesigns $d$ that can be obtained as in (a) or (b) above from the Clatworthy (1973) tables of PBIB designs, also accessible online via Google books. The results are summarized in Table 2, which again shows $k_0, k_{0.90}, k_{0.95}$ and $k_{0.98}$ for each such $d$. For each $(u, b)$ where multiple choices of $d$ emerge from the Clatworthy (1973) tables, we show in Table 2 only the one that comes with highest $A_{eff}(d_0)$, via (10), over a wide range of values of $k$. The details of the PBIB designs, such as R6, R10 etc., mentioned in Table 2 can be found from Clatworthy (1973).

For all the subdesigns $d$ in Table 2, except the last two, $A_{eff}(d_0)$ reaches 0.98 much before $k$ reaches

20. Of the two exceptional cases, the $d$ for $(u, b) = (18, 9)$ at least ensures $A_{\text{eff}}(d_0)$ reaching 0.95 quickly, while the $d$ for $(u, b) = (15, 10)$ has less impressive performance with $A_{\text{eff}}(d_0)$ increasing rather slowly beyond 0.90. Turning to the proportion, $f$, of singly replicated treatments, the growth of $f$ with $k$ is initially a little slow in Table 2 but catches up quickly. Among the 12 cases in this table, we have $f \geq 0.80$ in 4, 6 and 10 cases for $k = 15$, 20 and 25, respectively. We remark in this connection that even $k = 25$ can be relatively small compared to $v$, the total number of treatments, as $v$ also increases with $k$. Thus, for $k = 25$, we have $v \geq 80$ throughout Table 2, and $v \geq 126$ in 7 of the 12 cases there.

Table 2. *Values of $k_0$, $k_{0.90}$, $k_{0.95}$ and $k_{0.98}$ for subdesigns obtained from PBIB designs*

| $(u, b)$ | Subdesign $d$ | $k_0$ | $k_{0.90}$ | $k_{0.95}$ | $k_{0.98}$ | $(u, b)$ | Subdesign $d$ | $k_0$ | $k_{0.90}$ | $k_{0.95}$ | $k_{0.98}$ |
|---|---|---|---|---|---|---|---|---|---|---|---|
| (14, 4) | Dual of R6 | 8 | 8 | 8 | 8 | (21, 6) | Dual of R20 | 8 | 8 | 8 | 10 |
| (16, 4) | Dual of R10 | 9 | 9 | 9 | 9 | (24, 6) | Dual of R24 | 9 | 9 | 9 | 10 |
| (20, 4) | Dual of R16 | 11 | 11 | 11 | 11 | (16, 8) | LS28 or | | | | |
| (15, 5) | Dual of C6 | 7 | 7 | 7 | 8 | | dual of SR9 | 5 | 5 | 5 | 13 |
| (25, 5) | Dual of C9 | 11 | 11 | 11 | 11 | (18, 9) | Dual of LS1 | 5 | 5 | 8 | X |
| (12, 6) | Dual of R18 | 5 | 5 | 5 | 7 | (15, 10) | Dual of T2 | 4 | 6 | X | X |
| (18, 6) | Dual of R19 | 7 | 7 | 7 | 8 | | | | | | |

**Example 3**. Even for $u + b > 30$, use of PBIB designs can lead to high $A$-efficiency. As an illustration, let $(u, b) = (36, 12)$. Take the subdesign $d$ as LS74 in Clatworthy (1973). Then $(k_0, k_{0.90}, k_{0.95}, k_{0.98}) = (7, 7, 7, 11)$. This $d$, found painlessly as a known PBIB design, is isomorphic to the algorithmically obtained one in Haines (2021). We remark that our $A$-efficiency measure, that determines $(k_{0.90}, k_{0.95}, k_{0.98})$, is objectively based on the lower bound (9), as opposed to any empirical stopping rule in algorithmic construction, and makes no assumption about a constant block size for the competitors of $d$. In this example, $f$ equals 0.75 for $k = 15$, and 0.82 for $k = 20$. □

In the spirit of Example 2, the next example shows how "nearly" PBIB designs can work well as subdesigns $d$.

**Example 4**. (a) Let $(u, b) = (20, 6)$. Then $2u/b$ is not an integer and one cannot obtain $d$ directly as a PBIB design or its dual. We, therefore, start with the dual of the PBIB design R19 considered in Table 2. This dual has 18 treatments and 6 blocks as given by

$$\{1, 4, 7, 10, 13, 16\}, \{2, 5, 8, 11, 13, 17\}, \{3, 4, 9, 11, 14, 18\}$$
$$\{1, 6, 7, 12, 14, 17\}, \{2, 6, 8, 10, 15, 18\}, \{3, 5, 9, 12, 15, 16\}.$$

Add a treatment 19 to the first two blocks and another treatment 20 to the last two blocks of the above to get a subdesign $d$ in 20 treatments and 6 blocks. For the $d$ so obtained, we have $(k_0, k_{0.90}, k_{0.95}, k_{0.98}) = (7, 7, 7, 9)$, i.e., the associated $A_{\text{eff}}(d_0)$ reaches 0.98 very fast.

(b) We now explore a more severe modification of a PBIB design. Let $(u, b) = (20, 8)$. Then $2u/b$ is an

integer, but the Clatworthy (1973) tables do not show any PBIB design which or whose dual can be used as $d$ for the specified $(u, b)$. So, we start with the PBIB design LS28 mentioned in Table 2. It is a square lattice design with 16 treatments and 8 blocks as given by

$$\{1, 2, 3, 4\}, \{5, 6, 7, 8\}, \{9, 10, 11, 12\}, \{13, 14, 15, 16\},$$
$$\{1, 5, 9, 13\}, \{2, 6, 10, 14\}, \{3, 7, 11, 15\}, \{4, 8, 12, 16\}.$$

Suppose treatments 17, 18, 19 and 20 are added, respectively, to the (i) first two blocks, (ii) third and fourth blocks, (iii) fifth and sixth blocks, and (iv) last two blocks of the above. Then we get a subdesign $d$ in 20 treatments and 8 blocks. For this $d$, $(k_0, k_{0.90}, k_{0.95}, k_{0.98}) = (6, 6, 6, 10)$, i.e., the associated $A_{\text{eff}}(d_0)$ reaches 0.98 well before $k$ reaches 20.

In parts (a) and (b) above, $f$ equals 0.71 and 0.80 for $k = 15$, and 0.80 and 0.86 for $k = 20$. □

## 7. *MV*-efficiency

Consider any connected full design $d_0$. We continue to write $d$ for the subdesign of $d_0$ that contains the $u$ twice replicated treatments. As before, $W$ denotes the set of the $w$ singly replicated treatments, of which $s_1,\ldots, s_b$ appear in the $b$ blocks of $d_0$. Let

$$G = \{ j : s_j \geq 1; j = 1,\ldots, b\}. \tag{11}$$

By Proposition 1, then the *MV*-criterion for $d_0$ is given by

$$MV(d_0) = \max\{MV_{UU}, MV_{WW}, MV_{UW}\}, \tag{12}$$

where

$$MV_{UU} = \max\{(e_i - e_{i*})^T C^+ (e_i - e_{i*}): 1 \leq i < i* \leq u\}, \tag{13}$$

$$MV_{WW} = 2 + \max\{(\tilde{e}_j - \tilde{e}_{j*})^T \tilde{C}^+ (\tilde{e}_j - \tilde{e}_{j*})\}: j, j* \in G, j < j*\}, \tag{14}$$

$$MV_{UW} = (3/2) + \max\{\xi_{ij}^T \tilde{C}^+ \xi_{ij} : 1 \leq i \leq u, j \in G\}, \tag{15}$$

and $C^+$ and $\tilde{C}^+$ refer to the subdesign $d$ and its dual, respectively. We revisit Example 1 to illustrate (12)-(14) and also motivate our approach to studying the *MV*-criterion.

**Example 1** (continued). Here, $\tilde{C}^+ = \{2/(\lambda b)\}(I_b - b^{-1}J_{bb})$, $k \geq \lambda(b-1) + 1$, and $s_1,\ldots, s_b \geq 1$. Hence, $G = \{1,\ldots, b\}$, and by (14), $MV_{WW} = 2 + \{4/(\lambda b)\}$. Also, by (3), the entries in $\xi_{ij}$ sum to zero, and therefore,

$$(\lambda b/2) \xi_{ij}^T \tilde{C}^+ \xi_{ij} = \{\tilde{e}_j - (1/2)N^T e_i\}^T \{\tilde{e}_j - (1/2)N^T e_i\}$$

$$= 1 - e_i^T N \tilde{e}_j + (1/4)e_i^T NN^T e_i = 1 - n_{ij} + (1/2) \leq 3/2,$$

since $n_{ij} \geq 0$ and every diagonal element of $NN^T$ equals the common replication number, 2, in $d$. Thus,

by (15), $MV_{UW} \leq (3/2) + \{3/(\lambda b)\}$. Finally, by (13), $MV_{UU}$ does not exceed twice the largest eigenvalue of $C^+$, i.e., $MV_{UU} \leq 2(b-1)/b$. Therefore, in this example, $MV_{WW}$ is larger than both $MV_{UU}$ and $MV_{UW}$. So, by (12), $MV(d_0) = MV_{WW} = 2 + \{4/(\lambda b)\}$. □

The above outcome is by no means isolated. In many other situations too, $MV_{WW}$ exceeds both $MV_{UU}$ and $MV_{UW}$. This matches our intuition because $MV_{WW}$ is determined by the singly replicated treatments alone, and prompts us to derive lower bounds on $MV(d_0)$ via those on $MV_{WW}$ in Theorem 2 below. Some notation and a proposition will help. Given $u$, $b$ ($\geq 3$) and $k$, which fix $w = bk - 2u$, we continue to denote the class of connected full designs by $D_0 \equiv D_0(u, b, k)$. Let $k^\#$ be the smallest integer greater than or equal to $2u/(b-1)$, and write $k^+ = 2(u+1) - b$. Also, we continue to write $k_{\max}(d)$ for the largest block size in any subdesign $d$. Let $d_{MV}$ be $MV$-optimal among block designs with $b$ treatments and $u$ blocks, each of size two, and let $MV_{\min}$ denote the value of the $MV$-criterion for $d_{MV}$.

**Proposition 2**. (a) $k^+ > k^\#$. (b) *For every $d_0$ in $D_0$, with subdesign $d$, we have $k_{\max}(d) \leq k^+ - 1$.*

*Proof*. (a) Due to connectedness, $u \geq b - 1$. Hence (a) follows noting that

$$k^+ - \{2u/(b-1)\} = (b-2)[\{2u/(b-1)\} - 1] \geq 1.$$

(b) If $k_{\max}(d) \geq k^+$, then $d$ has at least $k^+$ plots in one block and its remaining $b - 1$ blocks together have at most $2u - k^+$ ($= b - 2$) plots, i.e., $d$ is not spread out over all $b$ blocks, which contradicts connectedness. This proves (b). □

**Theorem 2**. (a) *If $k \geq k^\#$, then for every $d_0$ in $D_0$, we have $MV(d_0) \geq 2 + 2\{(b-1)/u\}$.*

(b) *If $k \geq k^+$, then for every $d_0$ in $D_0$, we have $MV(d_0) \geq 2 + MV_{\min}$.*

*Proof*. Consider any $d_0$ in $D_0$ and the corresponding set $G$ in (11). Let $G$ have cardinality $g$, where $g \leq b$. Recall that because of connectedness,

$$s_j < k, \quad j = 1,\ldots, b. \tag{16}$$

(a) Suppose $k \geq k^\#$. Then $k \geq 2u/(b-1)$, so that

$$w = bk - 2u = \{k(b-1) - 2u\} + k \geq k. \tag{17}$$

Hence by (1), (11) and (16), $g \geq 2$. Without loss of generality, let $G = \{1, 2,\ldots, g\}$. Then $s_j = 0$ for $j > g$. So, by (1) and (16),

$$w = s_1 + \cdots + s_g < gk. \tag{18}$$

Now, let $I_g$ and $1_g$ denote, respectively, the identity matrix of order $g$ and the $g \times 1$ vector of ones, and define the nnd matrices

$$B = \begin{bmatrix} gI_g - 1_g 1_g^T & 0 \\ 0 & 0 \end{bmatrix}, \qquad B_1 = \begin{bmatrix} gI_g & 0 \\ 0 & 0 \end{bmatrix},$$

of order $b$. As $g \geq 2$, then by (14), as in the proof of Theorem 1,

$$MV_{WW} = 2 + \max_{1 \leq j < j^* \leq g} (\tilde{e}_j - \tilde{e}_{j^*})^T \tilde{C}^+ (\tilde{e}_j - \tilde{e}_{j^*})$$

$$\geq 2 + [2/\{g(g-1)\}] \Sigma\Sigma_{1 \leq j < j^* \leq g} (\tilde{e}_j - \tilde{e}_{j^*})^T \tilde{C}^+ (\tilde{e}_j - \tilde{e}_{j^*})$$

$$= 2 + [2/\{g(g-1)\}] \operatorname{tr}(\tilde{C}^+ B). \tag{19}$$

Next, observe that

$$\operatorname{tr}(B) - b^{-1} 1_b^T B 1_b = g(g-1),$$

while by (18) and Lemma 1(b),

$$\operatorname{tr}(\tilde{C}B) \leq \operatorname{tr}(\tilde{C}B_1) \leq (1/2) g \Sigma_{j=1}^g (k - s_j) = (1/2) g(gk - w).$$

Hence by Lemma 1(a),

$$(1/2) g(gk - w) \operatorname{tr}(\tilde{C}^+ B) \geq \{g(g-1)\}^2.$$

Because $gk - w > 0$ by (18), $w \geq k$ by (17), and $g \leq b$, $w = bk - 2u$, this yields

$$[2/\{g(g-1)\}] \operatorname{tr}(\tilde{C}^+ B) \geq 4(g-1)/(gk - w) \geq 4(b-1)/(bk - w) = 2(b-1)/u.$$

The truth of (a) is now immediate from (12) and (19).

(b) Suppose $k \geq k^+$. By Proposition 2(b), then $k_{\max}(d) < k$. Hence, by (11), $G = \{1,\ldots, b\}$ and by (14),

$$MV_{WW} = 2 + \max_{1 \leq j < j^* \leq b} (\tilde{e}_j - \tilde{e}_{j^*})^T \tilde{C}^+ (\tilde{e}_j - \tilde{e}_{j^*}) \geq 2 + MV_{\min},$$

which proves (b), in view of (12). The last inequality is evident from the definition of $MV_{\min}$, since the dual of $d$ has $b$ treatments and $u$ blocks, each of size two. □

Arguments similar to those in (19), in conjunction with the second inequality in (8), show that the second bound in Theorem 2 is at least as sharp as the first one there. Indeed, it is often sharper. On the other hand, the second bound can be used only when $MV_{\min}$ is known and $k \geq k^+$.

Along the lines of our approach to $A$-efficiency, we investigate how, for any choice of the subdesign $d$, Theorem 2 enables us to explore the $MV$-efficiency of the associated full design $d_0$ for various $k$. As in the proof of Theorem 2(b), if $k \geq k_{\max}(d) + 1$, then $G = \{1,\ldots, b\}$ and hence by (12)-(15), $MV(d_0)$ remains the same over $k$. Thus, $MV(d_0)$ stabilizes very fast with increase in $k$, since $k_{\max}(d) + 1$ either equals or is one more than the smallest possible block size of $d_0$, given $d$. Let $k^* [\equiv k^*(d)] = \max\{k^{\#}, k_{\max}(d) + 1\}$. In view of Theorem 2, the $MV$-efficiency of $d_0$, in the class $D_0$, can now be expressed as

$$MV_{\text{eff}}(d_0) = [2 + 2\{(b-1)/u\}]/MV(d_0), \quad \text{if } k \geq k^*, \tag{20}$$

$$= (2 + MV_{\min})/MV(d_0), \quad \text{if } k \geq k^+. \tag{21}$$

By Proposition 2, $k^+ \geq k^*$. Thus, both forms of $MV_{\text{eff}}(d_0)$ remain invariant over $k$ in the respective ranges, and this helps in the presentation of our findings.

**Example 1** (concluded). Here, $u = \lambda\binom{b}{2}$, $k^\# = \lambda b$, $k_{\max}(d) = \lambda(b-1)$, $k^* = \lambda b$ and $MV(d_0) = MV_{WW} = 2 + \{4/(\lambda b)\}$. Therefore, for $k \geq k^* (= \lambda b)$, by (20), $MV_{\text{eff}}(d_0) = 1$, i.e., $d_0$ is $MV$-optimal in $D_0$, in addition to enjoying high $A$-efficiency. □

Other subdesigns, including those in Tables 1, 2 and Examples 2-4, are revisited under the $MV$-criterion in Appendix 2. From the findings there, we reach the satisfying conclusion that a vast majority of these subdesigns entail full designs that have high $MV$-efficiencies or are even $MV$-optimal, in addition to having impressive performance under the $A$-criterion.

## 8. Two variants of the $A$-criterion

As seen at the end of the last section, by and large, there is no conflict between the $A$- and $MV$-criteria with regard to the choice of the subdesign $d$. A vast majority of the subdesigns studied here lead to high efficiencies under both criteria. Thus, a choice between the two criteria is not a source of major concern but essentially a matter of taste.

The findings in Section 7 and Appendix 2 show that the treatments in $W$ play a key role in determining the $MV$-criterion, which is natural because these are the ones that are singly replicated. From this perspective, we now focus on these treatments to consider two variants of the $A$-criterion based on the average variance of the BLUEs of (i) the $w(w-1)/2$ pairwise treatment contrasts arising from $W$ alone, and (ii) the $uw + w(w-1)/2$ pairwise treatment contrasts, arising from comparisons of the treatments in $W$ with the twice replicated treatments in $U$, in addition to those in (i). These are termed the $A_1$- and $A_2$-criteria. Clearly, a high efficiency under the $A_1$- or $A_2$-criterion calls for keeping, respectively, $A_{WW}$ or $A_{UW} + A_{WW}$ small. Here $A_{UW}$ and $A_{WW}$ are as in the proof of Theorem 1 and, following this proof,

$$A_{WW} = w(w-1) + \text{tr}\{\widetilde{C}^+(wS - ss^T)\},$$

$$A_{UW} + A_{WW} = (1/2)w(3u + 2w - b - 1) + (1/2)k\{w\text{tr}(\widetilde{C}^+) + b\text{tr}(\widetilde{C}^+ S)\},$$

in a full design $d_0$, associated with a subdesign $d$ to which $\widetilde{C}^+$ relates.

The term $ss^T$ appearing with a minus sign hinders finding a lower bound on $A_{WW}$ and necessitates the use of an enumerative approach for exploring the $A_1$-criterion. On the other hand, analogously to (9), a lower bound on $A_{UW} + A_{WW}$, over the class $D_0$ [$\equiv D_0(u, b, k)$] emerges as

$$A_{2\text{bound}} = (1/2)w(3u + 2w - b - 1) + kw\{(b-1)^2/u\},$$

leading to a measure of $A_2$-efficiency of any design $d_0$ in $D_0$ as $A_{2\text{eff}}(d_0) = A_{2\text{bound}}/(A_{UW} + A_{WW})$.

We are now in a position to reconsider the subdesigns in Sections 4-6, that is, those in Examples 1-4 and Tables 1 and 2, under the $A_1$- and $A_2$-criteria. As in Section 5, for $k = 10, 15$ and $20$, a complete enumeration shows that every $d_{\text{cat}}(u, b)$ in Table 1 minimizes both $A_{WW}$ and $A_{UW} + A_{WW}$ among binary competing subdesigns. Hence, the same holds in the entire competing class when $u = b - 1$. Such complete enumeration is difficult for the subdesigns $d$ in Examples 1-4 and Table 2. However, for $k \geq 10$, large scale random generation of competing designs yields a smaller $A_{WW}$ than these $d$ only in very rare situations, the reduction in $A_{WW}$ in these rare cases being typically less than 5%, and often less than 1%. Turning to the $A_2$-criterion, which can be handled analytically via the expressions for $A_{2\text{bound}}$ and $A_{2\text{eff}}(d_0)$ stated above, it is clear from Example 1 that if the subdesign $d$ is chosen as a linked block design, then the resulting full design $d_0$ is $A_2$-optimal. Furthermore, all $d$ in Examples 2-4 and Table 2 entail a little higher $A_{2\text{eff}}(d_0)$ than $A_{\text{eff}}(d_0)$ for each of the many values of $k$ that we examined.

The facts noted above are reassuring. They show that these subdesigns, that were earlier seen to perform well under the $A$-criterion, continue to do so under the two variants of this criterion considered here. The same pattern holds also for the large designs considered in the next section.

## 9. Large designs and a practical guideline

### 9.1 Large designs

It is of interest to see how far our results apply to early generation varietal trials where the number of treatments, $v$, can well exceed 1000. If $v$ is large, then so is $bk$, the total number of plots. Hence, if the number of twice replicated treatments, $u$, is too small, then the block size $k$ can be unrealistically large, as connectedness demands $u \geq b - 1$. Due to this reason, the designs from designtheory.org, discussed in Section 5, are not very useful for large $v$. In other words, if $v$ is large, then $u$ should also be commensurately large, e.g., if $v = 1500$, then one may consider $u = 100$ and still the proportion $f$ of singly replicated treatments remains as high as 0.93. As we shall see, the approaches in Sections 4 and 6, based on linked block and PBIB designs, can then be successfully employed. All large designs considered here satisfy $f \geq 0.91$ and $k/v \leq 0.09$. In fact, quite a few of these have appreciably larger $f$ or smaller $k/v$, e.g., the design in Table 3(g) below has $f = 0.96$ and $k/v = 0.06$.

For various $(u, b, k)$, each associated with a large $v = bk - u$, Table 3 illustrates the use of linked block and PBIB designs to obtain subdesigns $d$ such that the corresponding full designs $d_0$ have high $A$- and $MV$-efficiencies, as given by eff($d_0$) = ($A_{\text{eff}}(d_0)$, $MV_{\text{eff}}(d_0)$). The PBIB designs that give $d$ in (c)-(g) of this table can be found from Clatworthy (1973).

Table 3. *Subdesigns d leading to efficient large full designs $d_0$*

| ($u, b, k$) | $v$ | $d$ | eff($d_0$) |
|---|---|---|---|
| (a) (91, 14, 105) | 1379 | Linked block design | (0.999, 1) |
| (b) (190, 20, 120) | 2210 | Linked block design | (0.999, 1). |
| (c) (105, 21, 80) | 1575 | Dual of T8 | (0.993, 0.987) |
| (d) (100, 20, 80) | 1500 | LS 137 | (0.996, 0.992) |
| (e) (100, 25, 64) | 1500 | Dual of LS5 | (0.979, 0.969) |
| (f) (36, 12, 90) | 1044 | LS74 | (0.991, 0.979) |
| (g) (40, 16, 67) | 1032 | Dual of M1 | (0.971, 0.956) |

In Table 3, (d) and (e) give two possibilities for $d_0$ for the same ($u$, $v$) = (100, 1500) and hence the same $bk$ = 1600, with impressive eff($d_0$) in both cases. However, (e) achieves a 20% smaller $k$ than (d), whereas (d) ensures an even higher eff($d_0$) in its competing class than (e) does. Given $bk$ = 1600, if there is flexibility about $b$ and $k$, then a choice between (e) and (d) depends on whether or not one prioritizes a significantly smaller $k$ over somewhat a higher eff($d_0$). Elsewhere too, it may be possible to find alternative designs that reduce $k$ substantially, while retaining satisfactory eff($d_0$). As illustrated in (a*)-(c*) below, a type of extended group divisible (EGD) design, denoted by EGD($m_1,\ldots, m_p$) design and detailed in Appendix 3, can be conveniently used for this purpose. Any such EGD design is a PBIB design with more than two associate classes and hence beyond the range of the Clatworthy (1973) tables.

(a*) ($u, b, k$) = (91, 21, 70), $v$ = 1379: Augment EGD(7, 3) design by seven new blocks

{1, 5}, {4, 8}, {7 11}, {10 14}, {13 17}, {16 20}, {19 2},

and take $d$ as the dual of the resulting 91-block design; eff($d_0$) = (0.979, 0.956).

(b*) ($u, b, k$) = (190, 48, 50), $v$ = 2210: Delete the two blocks {1, 2} and {32, 48} from EGD(3, 4, 4) design and take $d$ as the dual of the resulting 190-block design; eff($d_0$) = (0.968, 0.933).

(c*) ($u, b, k$) = (105, 30, 56), $v$ = 1575: Take $d$ as the dual of EGD(2, 3, 5) design; eff($d_0$) = (0.961, 0.938).

The $d_0$ in (a*)-(c*) have high $A$-efficiencies. Their $MV$-efficiencies are also acceptable, especially as these relate to possibly non-attainable bounds. Like (e) and (d) in Table 3, the $d_0$ in (a*) has the same $u$, $v$ and $bk$ as that in Table 3(a), but with a 33% smaller $k$ and a lower but still high eff($d_0$). The picture is similar when one compares the $d_0$ in (b*) or (c*), respectively, with those in Table 3(b) or (c).

The idea of reducing block size with some sacrifice of efficiency can be carried out a step further if the numbers, $u$ and $w$, of twice and singly replicated treatments are not absolutely rigid but allow minor adjustments. For instance, the designs in (f) and (g) of Table 3 have similar ($u, w$), namely, (36, 1008) and (40, 992). Suppose one more replication is allowed for only four of the $w = 1008$ singly replicated treatments in (f), so that $(u, w)$ there becomes (40, 1004). On the other hand, if a new singly replicated treatment is added to each of the first 12 blocks of the $d_0$ in (g), then the resulting full design, say $d^*$, having block sizes 68 and 67, also has $(u, w) = (40, 1004)$. It can be a good substitute for the $d_0$ in (f), with block size reduced by about 25%, provided the minor adjustment in $u$ and $w$, as indicated above, is allowable in the latter. Although our expressions for efficiency do not hold for $d^*$ due to its unequal block sizes, one can intuitively expect $d^*$ to be highly efficient in the comparable class of designs, because of its close proximity to the $d_0$ in (g). Indeed, a comparison with appropriate bounds shows that $d^*$ has $A$-and $MV$-efficiencies 0.965 and 0.956 in the relevant class. The details, applicable only to certain block sizes given $(u, b)$, are omitted here to save space, but briefly touched upon later in Section 10.

As mentioned earlier, all the designs in this section have $f \geq 0.91$, i.e., the number, $u$, of twice replicated treatments is less than 10% of $v$, the total number of treatments. While such pre-specification of $u$ may be a practical necessity for experiments with large $v$, our approaches work equally well if the available resources allow a larger $u$ relative to $v$, leading to a smaller $f$. For example, if $(u, b, k) = (435, 30, 80)$, then $v = 1965$, $w = 1530$, $f = 0.78$, and with $d$ is chosen as a linked block design, the corresponding full design $d_0$ has eff$(d_0) = (0.999, 1)$. Similar examples abound.

*9.2 A practical guideline*

In this paper, we have presented several approaches for obtaining an efficient partially replicated full design $d_0$ by properly selecting its subdesign $d$ which consists of the twice replicated treatments. A guideline on how these approaches can be used systematically in a given context will be helpful. Given the numbers $u$ and $w$ of twice and singly replicated treatments, the number of blocks $b$, and the block size $k$ of $d_0$, which satisfy $bk = 2u + w$, the choice of $d$ depends on $u$ and $b$. On the basis of the results here, the possibilities in (i)-(iv) below may be explored for selecting $d$ with a view to obtaining an efficient $d_0$. Note that the provision of "modification" in (i), (iii) and (iv), as in Examples 2, 4 or (a*), (b*) of the previous subsection, considerably enhances their scope of applicability.

(i) $d$ chosen as a linked block design or a modification thereof (Subsection 4.1); applicable when $u$ equals or nearly equals $\lambda\binom{b}{2}$, for some positive integer $\lambda$.

(ii) $d$ chosen as the dual of an $A$-optimal block design with block size two as one can find from the website `designtheory.org` (Section 5); applicable when $u + b \leq 16$.

(iii) $d$ chosen as a two-replicate PBIB design or the dual of a PBIB design with block size two as one can find from the Clatworthy (1973) tables, or a modification thereof (Section 6); applicable when such a PBIB design is available from these tables, which can be quickly checked from Table 0 there.

(iv) $d$ chosen as the dual of an EGD design or a modification thereof (Subsection 9.1 and Appendix 3); applicable when $b = m_1 \ldots m_p$ and $u$ equals or nearly equals $b(m_1 + \ldots + m_p - p)/2$, for integers $p$ ($\geq 2$) and $m_1, \ldots, m_p$ ($\geq 2$).

As an illustration, let $(u, b, k) = (83, 27, 42)$, i.e., $w = 968$, which corresponds to a large $v = u + w = 1051$. Obviously, (i) or (ii) do not work here. Upon checking Table 0 of Clatworthy (1973), it is seen that neither a two-replicate PBIB design nor a PBIB design with block size two, coming close to the stipulated $(u, b)$ is available there. Thus (iii) is also not applicable. Turning to (iv), $b$ has only two possible factorizations, 9x3 and 3x3x3, giving $u = 135$ and 81, of which the second one comes close to the specified $u$. Therefore, we add two new blocks {1, 18} and {10, 27} to EGD(3, 3, 3) design and take $d$ as the dual of the 83-block design so obtained. For $k = 42$, the resulting $d_0$ has $A_{\text{eff}}(d_0) = 0.961$ and $MV_{\text{eff}}(d_0) = 0.933$. While the $A$-efficiency is quite high, the $MV$-efficiency may also be considered satisfactory because, as before, it relates to a possibly non-attainable bound. Here $f = 0.94$ and $k/v = 0.04$.

If none of (i)-(iv) works, then in the spirit of the last paragraph of the preceding subsection, one may examine if a minor adjustment of $u$ and $w$ is allowed and whether or not it helps. Finally, in the worst scenario when none of these theory driven approaches is applicable, one needs to consider an algorithmic construction, for instance, along the lines of Haines (2021) or using some package. There too, the bounds on the $A$- and $MV$-criteria obtained here, rather than any empirical stopping rule, may help, unless the average block size, $2u/b$, of $d$ is too small compared to $u$, when these bounds may become over-conservative. We further remark that algorithmic searches compete well with our approaches for relatively small $v = u + w$, but become increasingly slow with increase in $v$, and tend to break down for large $v$, like $v \geq 1000$, or $u > 100$. On the other hand, as seen in the previous subsection, even for large $v$, if any of (i), (iii) or (iv) works, then it yields an efficient full design effortlessly.

**10. Concluding remarks**

In this paper, we proposed several approaches, such as those based on linked block designs, an online design catalog, PBIB designs from Clatworthy (1973) and EGD designs, for obtaining efficient partially replicated block designs, under the $A$- and $MV$-criteria, when the singly replicated treatments far outnumber the twice replicated ones, due to resource constraints. These designs were also found to perform well under two variants of the $A$-criterion. Before concluding, we indicate some open problems that emerge from or are related to the present work.

We have not gone beyond block designs in this paper. This is because of the ease in implementation and resulting wide use of such designs; see Haines (2021) who emphasized this point, with a reference to Piepho and Williams (2016). However, we believe that some of our techniques, such as those for deriving the bounds on the *A*- and *MV*-criteria, may facilitate the study of efficient designs in more complex settings in the future. In particular, it will be of interest to examine if our findings can be linked with the row-column setup and mixed model in Cullis et al. (2006).

The framework of a constant block size for the full design $d_0$, as followed here, requires the number of plots, $2u + w$, to be an integral multiple of $b$. Else, the block sizes in $d_0$ have to be unequal. Although it should not be hard to extend Theorem 1 to this situation, the derivation of a design-independent lower bound on the *A*-criterion can be a challenge because some terms like $ss^T$ in its proof will no longer cancel out. This is akin to the difficulty encountered with $A_{WW}$ in Section 8. The bounds hinted at near the end of Subsection 9.1 in connection with such an unequal block full design apply only to a certain type of $u$, $b$ and block sizes, and development of more general results in this direction will be welcome. Even for full designs with equal block size, it will be of interest to find a sharp lower bound on $A_{WW}$.

Another area of research that deserves serious attention is the possible extension of our results to multiple environment trials where a different subset of treatments is replicated twice in each environment in such a way that each treatment is duplicated in the same number of environments. In this case, the same full design $d_0$ as obtained here may be considered for each environment, with a relabeling of the twice replicated treatments across environments. For example, with $(u, b, k) = (20, 6, 20)$, i.e., $w = 80$, $v = 100$, and five environments $E_1, \ldots, E_5$, one may wish to proceed as in Example 4(a) to get the common $d_0$ for all environments, and then label the twice replicated treatments as $1, \ldots, 20$ in $E_1$, $21, \ldots, 40$ in $E_2$, and so on. While this appears to be prospective, a rigorous investigation even under a fixed effects model will require accounting for underlying details such as the environment effects and possible treatment versus environment interaction. Further complexities may arise if certain effects are random.

We conclude with the hope that the present endeavor will generate interest in the above and related problems.

**Appendix 1: More on *A*-efficient designs**

(A) Bailey and Cameron (2013) obtained *A*-optimal block designs in blocks of size two when the number of blocks is equal to or one less than the number of treatments. Some of these designs appear in the catalog `v-b-k-phi`. As seen in Section 5 by complete enumeration, their duals, namely, $d_{\text{cat}}(b, b)$ ($b = 4, \ldots, 8$) and $d_{\text{cat}}(b - 1, b)$ ($b = 5, \ldots, 8$), entail *A*-optimality of $d_0$ at least among binary competitors, even when their $A_{\text{eff}}(d_0)$ values are not impressive. The same is very likely to be the case for the other *A*-optimal designs in Bailey and Cameron (2013). For instance, consider their *A*-optimal design for 9 treatments in 8 blocks each of size two. The subdesign $d(8, 9)$ obtained as its dual has 8 treatments and 9 blocks, and leads to $A_{\text{eff}}(d_0)$ only 0.719 for $k = 20$, which is not surprising because as many as 8 blocks of $d(8, 9)$ have size one each. Here $u = b - 1$, all connected designs have to be binary, and an enumeration of all of them again shows that $d(8, 9)$ actually gives an $A$-optimal full design $d_0$ for each of the many values of $k$ that we examined. While such complete enumeration is difficult for the designs in Bailey and Cameron (2013) which involve 9 or more blocks, large scale random generation of competing designs reveals no improvement over their dual designs under the $A$-criterion.

(B) The catalog `v-b-k-phi` includes six $A$-optimal block designs, namely, 6-4-3, 8-4-4, 9-6-3, 10-4-5, 10-5-4 and 12-4-6, that are two-replicate and hence themselves qualify as the subdesign $d$. However, these are identical to $d_{\text{cat}}(u, b)$ for $(u, b) = (6, 4), (8, 4), (9, 6), (10, 4), (10, 5)$ and $(12, 4)$, respectively, as considered in Table 1. Hence there is no need to study them separately.

(C) In continuation of Section 6, while looking for efficient subdesigns $d$ with $u + b > 16$, one may also consider $E$-optimal block designs cataloged in `designtheory.org` over a range wider than the listing of $A$-optimal block designs there. Since $E$-optimal designs maximize the smallest eigenvalue of the intrablock information matrix, they often perform well under the $A$-criterion too. Thus, the two-replicate ones among them merit consideration as the subdesign $d$. An examination of all such $E$-optimal designs in `designtheory.org`, which can be used as $d$ for $u, b \geq 4$, $u + b > 16$, shows that only two of them are not covered by Example 1 or Table 2. These correspond to $(u, b) = (12, 8)$ and $(14, 7)$. The first one has $(k_0, k_{0.90}, k_{0.95}, k_{0.98}) = (4, 5, X, X)$. The second one is better, with $(k_0, k_{0.90}, k_{0.95}, k_{0.98}) = (5, 6, 6, X)$.

(D) Figure 1(a) of Cullis et al. (2006) gives the layout of a design in 5 columns by 30 rows, involving $u = 30$ twice and $w = 90$ singly replicated treatments. The twice replicated treatments are not spread equally among either the columns or the rows. If the columns of this figure are interpreted as blocks, then the full design has $b = 5$ blocks each of size $k = 30$, while the subdesign, $d$, formed by the $u = 30$ twice replicated treatments, has $b = 5$ blocks of sizes 12, 12, 13, 12 and 11. For efficient construction of such a subdesign $d$, one may start with a linked block design as in Example 1, with $\lambda = 3$ and $b = 5$, identify any treatment therein that appears in the fifth but not in the third block, and then shift this treatment from the fifth to the third block. The resulting $d$ leads to a full design $d_0$ with $A_{\text{eff}}(d_0) = 0.9983$, for $k = 30$. On the other hand, if the rows of Figure 1(a) of Cullis et al. (2006) are taken as blocks, then the resulting subdesign has average block size two and two blocks of size only one. Hence, as with some designs in Table 1, our $A$-efficiency measure then becomes over-conservative.

# Appendix 2: *MV*-efficiency of *A*-efficient designs

We revisit subdesigns, including those in Tables 1, 2 and Examples 2-4, under the *MV*-criterion. Computations show that, as with Example 1, they all satisfy $MV(d_0) = MV_{WW}$ for $k \geq k^*$. As a result, by (14) and the fact that $G = \{1,\ldots, b\}$ whenever $k \geq k^*$, for every such subdesign $d$, we have

$$MV(d_0) = 2 + MV_{\text{dual}}, \tag{A.1}$$

for $k \geq k^*$, and hence for $k \geq k^+$, where $MV_{\text{dual}}$ is the value of the $MV$-criterion for the dual $\tilde{d}$. Equation (A.1) will be invoked while applying (20), (21) or Theorem 2 to these subdesigns.

(A) Consider now the subdesigns $d_{\text{cat}}(u, b)$ in Table 1. The catalog v-b-k-MV in designtheory.org facilitates their evaluation under the $MV$-criterion. This catalog shows $MV$-optimal block designs in the binary class, and the ones among them with block size two are $MV$-optimal also in the general class. The subdesigns $d_{\text{cat}}(u, b)$, except those for

$$(u, b) = (5, 4), (9, 4), (11, 4), (9, 5), (8, 6), (9, 7), \tag{A.2}$$

have duals, all with block size two, that appear in the catalog v-b-k-MV and hence are $MV$-optimal block designs. Thus, by (A.1), for $k \geq k^+$, the corresponding full designs $d_0$ attain the lower bound in Theorem 2(b) and are $MV$-optimal in $D_0$. As for the six exceptional cases in (A.2), it can be seen from (21) and (A.1) that the $d_{\text{cat}}(u, b)$ there lead to $MV_{\text{eff}}(d_0) = 0.929, 0.989, 0.984, 0.956, 0.889$ and $0.934$, for $k \geq k^+$. In fact, use of (20) shows that some of the $d_{\text{cat}}(u, b)$ ensure high $MV$-efficiency or even $MV$-optimality of $d_0$ for still smaller $k$ as well. For instance, the $d_{\text{cat}}(u, b)$ for $(u, b) = (8, 4), (10, 4), (11, 4)$ and $(11, 5)$ ensure $MV_{\text{eff}}(d_0) = 0.971, 0.975, 0.954$ and $0.974$ for $k \geq 6, 7, 8$ and $6$, respectively, whereas the $d_{\text{cat}}(u, b)$ for $(u, b) = (6, 4), (12, 4)$ and $(10, 5)$, also covered by Example 1, entail $MV$-optimality of $d_0$ for $k \geq 4, 8$ and $5$. We further remark that for each of the six $(u, b)$ in (A.2), if instead of $d_{\text{cat}}(u, b)$, the dual of the $MV$-optimal block design in v-b-k-MV is taken as the subdesign $d$, then by (A.1) and Theorem 2(b), the resulting $d_0$ is $MV$-optimal in $D_0$ for $k \geq k^+$, but this comes at the cost of some loss in $A$-efficiency.

(B) Turning to the subdesigns in Table 2, we show in Table A.1 the corresponding $MV_{\text{eff}}(d_0)$, for $k \geq k^*$, as obtained from (20). Most of these $MV$-efficiencies are very high. Furthermore, by Theorem 2(b) along with (A.1), the subdesigns for $(u, b) = (14, 4), (16, 4), (20, 4), (15, 5)$ and $(12, 6)$ actually lead to $MV$-optimality of $d_0$ for $k \geq k^+$, because their duals are $MV$-optimal block designs, as can be verified by enumerating all competing designs. For the other $(u, b)$ in Table 3, such complete enumeration is hard, but we anticipate the same outcome for many of them.

Table A.1. *Values of $k^*$ and $MV_{\text{eff}}(d_0), k \geq k^*$, for subdesigns obtained from PBIB designs*

| $(u, b)$ | $k^*$ | $MV_{\text{eff}}(d_0)$ | $(u, b)$ | $k^*$ | $MV_{\text{eff}}(d_0)$ | $(u, b)$ | $k^*$ | $MV_{\text{eff}}(d_0)$ |
|---|---|---|---|---|---|---|---|---|
| (14, 4) | 10 | 0.991 | (25, 5) | 13 | 0.993 | (24, 6) | 10 | 0.989 |
| (16, 4) | 11 | 0.990 | (12, 6) | 5 | 0.944 | (16, 8) | 5 | 0.958 |
| (20, 4) | 14 | 0.996 | (18, 6) | 8 | 0.989 | (18, 9) | 5 | 0.929 |
| (15, 5) | 8 | 0.981 | (21, 6) | 9 | 0.983 | (15, 10) | 4 | 0.889 |

(C) We finally consider the subdesigns $d$ in Examples 2-4. The $MV$-efficiencies of the associated full designs $d_0$ are obtained via (20) and shown below.

(i) $(u, b) = (13, 6)$: $d$ is as in Example 2(a); $MV_{eff}(d_0) = 0.923$, for $k \geq 6$.

(ii) $(u, b) = (17, 6)$: $d$ is as in Example 2(b); $MV_{eff}(d_0) = 0.971$, for $k \geq 7$.

(iii) $(u, b) = (36, 12)$: $d$ is as in Example 3; $MV_{eff}(d_0) = 0.979$, for $k \geq 7$.

(iv) $(u, b) = (20, 6)$: $d$ is as in Example 4(a); $MV_{eff}(d_0) = 0.968$, for $k \geq 8$.

(v) $(u, b) = (20, 8)$: $d$ is as in Example 4(b); $MV_{eff}(d_0) = 0.953$, for $k \geq 6$.

As noted earlier, the subdesign in Example 2(c) is the same as $d_{cat}(8, 4)$, which has already been covered in the context of Table 1.

**Appendix 3: One type of EGD design**

An EGD$(m_1,\ldots, m_p)$ design, where $p$ ($\geq 2$) and $m_1,\ldots, m_p$ ($\geq 2$) are integers, involves (I) $b = m_1 \ldots m_p$ treatments $1,\ldots, b$, arranged in the natural order as an $m_1 \times \ldots \times m_p$ hyper-rectangle $H$, and (II) $u = b(m_1 + \ldots + m_p - p)/2$ blocks, each of size two, as given by all ordered pairs arising from each one-dimensional slice of $H$. To be more specific, the arrangement in the natural order in (I) involves indexing the cells of $H$ by $p$-tuples $(h_1,\ldots, h_p)$, where $h_j = 1,\ldots, m_j$ for $j = 1,\ldots, p$, and then assigning treatments $1,\ldots, b$ to the lexicographically ordered cells. Also, the one-dimensional slices in (II) are like the rows and columns of a rectangle.

For example, if $p = 3$ and $m_1 = m_2 = 2$, $m_3 = 3$, then $b = 12$ and treatments $1, \ldots, 12$ are assigned, respectively, to the lexicographically ordered cells

(1,1,1), (1,1,2), (1,1,3), (1,2,1), (1,2,2), (1,2,3), (2,1,1), (2,1,2), (2,1,3), (2,2,1), (2,2,2), (2,2,3)

of $H$, leading to the arrangement

Layer 1:  1  2  3 ,    Layer 2:   7  8  9 .
          4  5  6                10 11 12

The one-dimensional slices are then (i) $\{1 + 3a, 2 + 3a, 3 + 3a\}$, $a = 0, 1, 2, 3$, (ii) $\{1 + a, 4 + a\}$, $a = 0, 1, 2, 6, 7, 8$, and (iii) $\{1 + a, 7 + a\}$, $a = 0, 1, \ldots, 5$. The three ordered pairs from each slice in (i) form three blocks, while each slice in (ii) and (iii) yields only one ordered pair and hence only one block. Thus, we have 24 [$= b(m_1 + \ldots + m_p - p)/2$] blocks altogether. The four EGD designs in Section 9 are larger but can be constructed in the same way; in fact, for one of these, $H$ is simply a 7×3 rectangle.

We remark that an EGD$(m_1,\ldots, m_p)$ design, as considered here, is a PBIB design based on a factorial association scheme, see Bailey (2004, pp. 344-347) for more details.

**Acknowledgement**: The work of RM was supported by a grant from Anusandhan National Research Foundation, Government of India, under their National Science Chair scheme.